\documentclass[a4paper,12pt]{article}

\usepackage{amsmath,amssymb}

\title{Darboux-Halphen-Brioshi system with rank four}
\date{}

\author{Y.~{Ohyama} and S.~{Okumura}\\
Department of Mathematics, Osaka University} 

\begin{document} 
\maketitle
\begin{abstract} 
We will determine the condition that a homogeneous quadratic equation 
with rank four becomes the Darboux-Halphen-Brioshi system in terms of 
the associated non-associative algebra. 
\end{abstract}

\def\bp{\mathbb{P}}
\def\bc{\mathbb{C}}
\def\bz{\mathbb{Z}}
\def\ds{\displaystyle}

\newtheorem{theorem}{Theorem}[section]
\newtheorem{proposition}[theorem]{Proposition}
\newtheorem{lemma}[theorem]{Lemma}
\newtheorem{corollary}[theorem]{Corollary}
\newtheorem{definition}[theorem]{%
\fontfamily{cmr}\fontseries{m}\fontshape{sc}\selectfont%
Definition}
\newtheorem{remark}[theorem]{%
\fontfamily{cmr}\fontseries{m}\fontshape{sc}\selectfont%
Remark}
\newtheorem{example}[theorem]{%
\fontfamily{cmr}\fontseries{m}\fontshape{sc}\selectfont%
Example}
\newenvironment{proof}{\noindent\textrm{Proof.}\\%
\indent}{\vspace{1ex}}

\hyphenation{qua-drat-ic hyp-er-geo-met-ric ho-mo-ge-ne-ous}


\section{Introduction}
We will study a homogeneous quadratic nonlinear differential equation
\begin{equation}\label{quad}
\frac{d X_i}{dt}=\sum_{j,k=1}^n a_{jk}^iX_j X_k.
\end{equation} 
This type of equations contains many important examples such as 
Euler's spinning top motion equation or the Lotka-Volterra equation (see 
Example \ref{Euler-top} and \ref{Lotka-Volterra}). 

In recent years, another class of quadratic equations called Halphen's 
first equation \cite{H1} 
\begin{align}\label{hal1}
  \frac{dX_2}{dt}+\frac{dX_3}{dt} &= 2 X_2 X_3,\nonumber\\
  \frac{dX_1}{dt}+\frac{dX_3}{dt} &= 2 X_1 X_3,\\
  \frac{dX_1}{dt}+\frac{dX_2}{dt} &= 2 X_1 X_2 \nonumber
\end{align}
is studied in many areas, such as the Atiyah-Hitchin metrics, modular 
forms, moonshine \cite{HM}, 
or the Painlev\'e analysis.  Halphen's first equation is satisfied by 
the null values of elliptic theta functions \cite{O}. Since the null 
values of elliptic theta functions are related to the hypergeometric 
equations, Halphen also studied Halphen's second equation in 
1881~\cite{H2}: 
\begin{align}\label{hal2}
\frac{dX}{d\tau} &= X^2 +c (X-Y)^2 +b (Z-X)^2+a(Y-Z)^2,\nonumber\\
\frac{dY}{d\tau} &= Y^2 +c (X-Y)^2 +b (Z-X)^2+a(Y-Z)^2,\\
\frac{dZ}{d\tau} &= Z^2 +c (X-Y)^2 +b (Z-X)^2+a(Y-Z)^2,\nonumber
\end{align}
which is solved by hypergeometric functions. 
Halphen's second equation is the simplest example of the 
Darboux-Halphen-Brioshi system \cite{O1} (see the section 4). 

In this paper we will study the condition when quadratic equations are 
reduced to the Darboux-Halphen-Brioshi system. We will write this 
condition in terms of non-associative algebras, which are found by 
L.~Markus \cite{M}. 

In the case of rank three, this condition is completely determined in 
\cite{O2}. In this case the corresponding Darboux-Halphen-Brioshi system 
is related to hypergeometric equations or confluent hypergeometric 
equations. 

The rank four case is important in the study of Painlev\'e equations, 
since the Painlev\'e sixth equations represent isomonodromic 
deformations of linear differential equations with four regular 
singularities. Moreover one of the simplest Darboux-Halphen-Brioshi 
system with rank four is related to modular forms of level three 
\cite{O3}. This system is also related to the quantum cohomology of 
$\bc\bp^2$. 

The first author (Y.~O.) thanks to Professor John Harnad and Professor 
John McKay for fruitful discussions.

\section{Homogeneous Quadratic Differential Equations}
We take a homogeneous quadratic nonlinear differential equations 
\begin{equation}
\frac{d X_i}{dt}=\sum_{j,k=1}^n a_{jk}^iX_j X_k,
\end{equation}
where $a_{jk}=a_{kj}$.
This type of equations contains many important examples. 

\begin{example} \label{Euler-top}
  Euler's spinning top motion equation:
  \begin{align*}
  \frac{dX_1}{dt}&= 2 X_2 X_3,\\
  \frac{dX_2}{dt}&= 2 X_1 X_3,\\
  \frac{dX_3}{dt}&= 2 X_1 X_2.
  \end{align*}
\end{example}

\begin{remark}
Even if a quadratic equation is not 
homogeneous, we can interpret it into homogeneous equation with 
dummy variables. 
\end{remark}

\begin{example} \label{Lotka-Volterra}
  The Lotka-Volterra equation
  \begin{align*}
  \frac{dN_1}{dt}&= a N_1-b N_1 N_2,\\
  \frac{dN_2}{dt}&= -c N_1-d N_1 N_2,
  \end{align*}
is an inhomogeneous equation, but, with a dummy variable $N_3$, it can 
be represented as a homogeneous equation:
 \begin{align*}
  \frac{dN_1}{dt}&= a N_1N_3-b N_1N_2,\\
  \frac{dN_2}{dt}&= -c N_1N_3-d N_1N_2,\\
 \frac{dN_3}{dt}&=0.
 \end{align*}
\end{example}


We define a non-associative algebra ${\cal A}$ 
\begin{equation*}
{\cal A}=\sum_{k=1}^n \bc \ x_k, \qquad x_j \cdot x_k =\sum_{i=1}^n 
a_{jk}^i \ x_i ,
\end{equation*}
which is associated with the homogeneous equation \eqref{quad}. Since 
$a_{jk}=a_{kj}$, ${\cal A}$ is commutative.
\par\medskip

\begin{remark}
A relation between non-associative algebra and homogeneous equation is 
found by L.~Markus in 1960 \cite{M}. That is used in the study of 
classification of topological behavior of solution orbits near the 
origin, which is a critical point of \eqref{quad} . 
\end{remark}

From now on, we use uppercase letters for independent variables of a 
homogeneous quadratic equation, and lowercase letters for basis of the 
corresponding commutative non-associative algebra. 
The independent variables are dual to bases of the corresponding algebra: 
\begin{lemma}
The commutative non-associative algebra ${\cal A}$ and the homogeneous 
quadratic nonlinear differential equation admits a one to one 
correspondence. 
Moreover, ${\cal A}$ determines a structure of homogeneous quadratic 
equation. If and only if two homogeneous quadratic nonlinear 
differential equations 
\begin{align*}
\frac{d X_i}{ dt}&=\sum_{j,k=1}^n a_{jk}^iX_j X_k,&
\frac{d Y_i}{ dt}&=\sum_{j,k=1}^n b_{jk}^iY_j Y_k,
\end{align*}
are interchanged by a linear transformation
\begin{equation*}
Y_j =\sum_{k=1}^n c_{jk} X_k,
\end{equation*}
then the corresponding commutative non-associative algebras are 
equivalent: 
\begin{equation*}
x_k = \sum_{k=1}^n c_{jk} y_j.
\end{equation*}
\end{lemma}

The lemma above is easily verified.

\section{Halphen's Equation and Hypergeometric Function}

The Halphen's second equation \eqref{hal2} is solved by a 
hypergeometric functions.

\begin{remark}
If $a=b=c=-1/8$, then the equation \eqref{hal2} turns into the Halphen's 
first equation \eqref{hal1}, and then it is solved by theta functions.
\end{remark}

In \cite{O2} we demonstrated the condition that a third order quadratic 
system is solved by hypergeometric functions. 

\begin{theorem} \label{rank3}
If and only if a three dimensional commutative non-associative algebra admits the 
unit, then the corresponding quadratic system is solved by (confluent) 
hypergeometric functions or elementary functions.  

If the automorphism group of the non-associative algebra is a finite 
group, the quadratic equation is solved by (confluent) 
hypergeometric functions. If the automorphism group of the 
non-associative algebra is an infinite group, the quadratic equation 
is solved by elementary functions. 
\end{theorem}

\begin{example}[Matrix Riccati equation]
Let $A$, $X$ be $2\times 2$ symmetric matrices which 
satisfy 
\[
\frac{dX}{dt}=XAX. 
\]
If $A$ is a constant matrix, then the equation is a third order 
quadratic system, and then the associated commutative non-associative 
algebra is a Jordan algebra. Since the algebra admits the unit and the 
automorphism is $O(2)$, the equation is solved by elementary functions. 
This algebra is independent of $A$ upto an isomorphism, if  $\det A \neq 0$. 
\end{example}

\begin{example}
For the Halphen's second equation \eqref{hal2}, we take the basis of the 
associated algebra as 
\begin{align*}
 e_1&= -x+y+z, &e_2&= x-y+z, &e_3&= x+y-z.
\end{align*}
Then the multiplication is as follows.
\begin{align*} 
e_1\cdot e_1 &= (1+4(b+c))\,e, &e_1\cdot e_2 &= -e_3-4\,c\, e,\\
e_2\cdot e_2 &= (1+4(a+c))\,e, &e_2\cdot e_3 &= -e_1-4\,a\, e,\\
e_3\cdot e_3 &= (1+4(a+b))\,e, &e_1\cdot e_3 &= -e_2-4\,b\, e, 
\end{align*}
where $e=e_1+e_2+e_3$. This algebra has the unit $e$. 
Then \eqref{hal2} is solved by the hypergeometric function $F(\alpha, 
\beta, \gamma, z)$, where
\begin{align*} 
a &= \frac14(2\alpha\beta-\gamma-\alpha\gamma-\beta\gamma+\gamma^2),\\
b &= \frac14(\alpha^2+\beta^2+\gamma-\alpha\gamma-\beta\gamma-1),\\
c &= \frac14(-2\alpha\beta-\gamma+\alpha\gamma+\beta\gamma).
\end{align*}
Let $p,q,r$ be an exponents of hypergeometric equation at singular 
points
\begin{align*}
 p&=1-\gamma, &q&= -\alpha -\beta + \gamma, &r&=\alpha -\beta, 
\end{align*}
then 
\begin{align*}
 1+4(a+c)&=p^2, &1+4(a+b)&=q^2, &1+4(b+c)&=r^2.
\end{align*}
\end{example}
\par\bigskip

If 
\begin{align*}
\alpha&=\frac{k+9}{12k}, &\beta&=\frac{k-9}{12k}, &\gamma&=\frac12,
\end{align*}
then by the transformation given by Chazy \cite{C}, Halphen's second 
equation \eqref{hal2} turns into 
\begin{equation}\label{chazy}
y^{\prime\prime\prime}=2y^{\prime\prime}-3(y^{\prime})^2+
\frac4{36-k^2}(6y^{\prime}-y^2)^2. 
\end{equation}
Moreover, if $k=0$, then the equation \eqref{chazy} is equivalent with 
\begin{align*} 
\frac{dX}{dt} &= X^2 +(V-X)(W-X), \\
\frac{dW}{dt} &= W^2-(X-W)^2+(V-X)(W-X), \\
\frac{dV}{dt} &=V^2+(W-X)^2-(X-V)^2+(V-X)(W-X), 
\end{align*}
which is solved by Airy functions \cite{O1}. 
This is also commented in \cite{CO}.

\section{Generalized Darboux-Halphen-Brioshi System}
For a Fuchsian equation
$$ \frac{d^2y}{ dz^2} + Q(z)y=0,$$
$$
Q(z)=\sum_{j=1}^m \frac{\alpha_j }{ \left(z-a_j\right)^2}+ 
\sum_{j=1}^{m-1} \frac{\beta_j}{ 
\left(z-a_j\right)\left(z-a_{j+1}\right)}, 
$$
if we set 
\begin{align*}
\tau &=\frac{y_2}{y_1}, & X_0 &= \frac{d}{d\tau}\log y_1, & X_j = 
\frac{d}{d\tau}\log \frac{y_1}{x-a_j},
\end{align*}
then
\begin{gather}\label{gDHB}
\frac{d{X_k}}{ d\tau} = X_k^2 - \sum_{j=1}^m \alpha_j \left( 
X_j-X_0\right)^2 
-\sum_{j=1}^{m-1} \beta_j 
\left(X_j-X_0\right)\left(X_{j+1}-X_0\right),\\ 
(X_j, X_k, X_l, X_n)=(a_j, a_k, a_l, a_n), \nonumber 
\end{gather}
where $(a, b, c, d)$ is an anharmonic ratio
$$(a, b, c, d)=\frac{a-b}{c-d}\frac{c-b}{a-d}.$$
\par\bigskip

We call \eqref{gDHB} the generalized Darboux-Halphen-Brioshi system 
(together with the algebraic relations) \cite{O2}. 
The variable $X_j$ is called the Brioshi variable.

\begin{remark}
If $m=2$, there is no algebraic relation. In this case the generalized 
Darboux-Halphen-Brioshi system is the same as Halphen's second equation. 
\end{remark}

The generalized Darboux-Halphen-Brioshi system was found in the study of 
differential relations of modular forms \cite{O3}. It plays 
an important role in the study of moonshine \cite{HM}. 

\begin{example}[The level-three Halphen equation] 
The equation
$$
\left\{\aligned
W^{\prime}+X^{\prime}+Y^{\prime}&=WX+XY+YW,\\
W^{\prime}+Y^{\prime}+Z^{\prime}&=WY+YZ+ZW,\\
W^{\prime}+X^{\prime}+Z^{\prime}&=WX+XZ+ZW,\\
X^{\prime}+Y^{\prime}+Z^{\prime}&=XY+YZ+ZX,\\
e^{\frac 4 3\pi i}(XZ+YW)&+e^{\frac 2 3\pi i}(XW+YZ)+(XY+ZW)=0
\endaligned \right. $$
is called the level-three Halphen equation. 
This is a generalized Darboux-Halphen-Brioshi system given by
the Picard-Fuchs equation 
\begin{equation}\label{pf-l3}
(1-t^3)y^{\prime\prime}-3t^2 y^{\prime}-ty=0
\end{equation}
of the $\Gamma(3)$-modular surface
$$x^3+y^3+z^3-3 t xyz=0.$$
The equation \eqref{pf-l3} is the Laplace transform of a linear equation 
$$\left( z\frac{d}{dz}\right)^3\phi =27 z^3\phi,$$
which appear in the Frobenius structure \cite{Dub} which turns up in the 
quantum cohomology of $\mathbb{CP}^2$. 
\end{example}

Now, we will study a quadric homogeneous equation with rank four
\begin{equation*}({\cal D}_3)\ 
\begin{cases}
\ds \frac{d{X_i}}{d\tau} = \sum_{j,k=0}^3 a_{jk}^iX_j X_k,\\
\ds Q(X_0, X_1, X_2, X_3)= \sum_{j,k=0}^3 b_{jk}X_j X_k,
\end{cases}
\end{equation*}
where $a_{jk}=a_{kj}, b_{jk}=b_{kj}$. 
We assume a compatibility condition
$$\frac{d}{d\tau}Q(X_0, X_1, X_2, X_3) = 
L(X_0, X_1, X_2, X_3)\,Q(X_0, X_1, X_2, X_3).$$
This condition means that the hypersurface $Q=0$ is invariant 
under the vector field defined by $({\cal D}_3)$.
If $L=0$ then $Q$ is a first integral. From now on, we assume $L\neq 0$.
Moreover, we assume that $Q$ is irreducible. 
\par\bigskip

We define a commutative non-associative algebra ${\cal A}_3(c)$ with a 
parameter $c=(c_1,c_2,c_3,c_4)$, which is associated with ${\cal D}_3$, 
as 
\begin{gather*}
{\cal A}_3(c)=\sum_{j=0}^3 \bc x_j,\\
 x_j\cdot x_k= \sum_{i=0}^3 a_{jk}^i x_i + b_{jk} x, \qquad x= 
 \sum_{i=0}^3 c_i x_i. 
\end{gather*}
We use this algebra ${\cal A}_3(c)$ to study the condition that the 
system ${\cal D}_3$ is a generalized Darboux-Halphen-Brioshi system.

\begin{theorem}
The system ${\cal D}_3$ is a generalized Darboux-Halphen-Brioshi system, 
if and only if there exist basis $x_0,x_1,x_2,x_3$ of ${\cal A}_3(c)$ 
which satisfies the following conditions: 
\begin{enumerate}
\item $e=x_0+x_1+x_2+x_3$ is the unit 
for any $c$. 
\item $(\pm x_0\pm x_1\pm x_2\pm x_3)^2$ is proportional to the unit $e$ for 
some $c$. 
\item $(- x_0 +x_1 +x_2 +x_3)^2$, $( x_0 -x_1 +x_2 +x_3)^2$, 
 $( x_0 +x_1 -x_2 +x_3)^2$, $( x_0 +x_1 +x_2 -x_3)^2$ are proportional 
 to the unit $e$ for any $c$.
\end{enumerate}
\end{theorem}

In the light of the algebra associated with Halphen's second equation, 
this theorem is an extension of theorem \ref{rank3}.

\par\bigskip
\begin{proof}
We set a basis of the algebra ${\cal A}_3(c)$ associated with the 
system ${\cal D}_3$: 
\begin{align*} 
e_0 &= x_0+x_1+x_2+x_3,\\
e_1 &= x_0-x_1+x_2+x_3,\\
e_2 &= x_0+x_1-x_2+x_3,\\
e_3 &= x_0+x_1+x_2-x_3.
\end{align*}
Since $e=e_0$ is the unit for any $c$, we obtain $Q=Q(E_1,E_2,E_3)$. 
Moreover, $e_j^2$ and $ (x_0-x_1-x_2-x_3)^2$ are proportional to the 
unit $e$ for any $c$, therefore we obtain 
$$Q=\gamma_1E_1 E_2+\gamma_2 E_2 E_3+ \gamma_3 E_1 E_3,$$
where $\gamma_1+\gamma_2+\gamma_3=0$. 
Since $Q$ is irreducible, we obtain $\gamma_j\neq 0\ (j=1,2,3)$. 
We set 
\begin{equation}\label{ee1}
e_j^2= \widetilde{\alpha_j}\,e_0\quad (j=1,2,3).
\end{equation}
Then, since we have
$$( x_0+ x_1-x_2- x_3)^2=( e_2+ e_3 - e_0)^2 \equiv
 e_2\cdot e_3 -(e_2+ e_3) \pmod{e_0},$$
we obtain 
$$ e_2\cdot e_3 \equiv e_2+ e_3 \pmod{e_0}.$$
In the same way, we obtain
\begin{align*}
e_1\cdot e_2 \equiv e_1+ e_2 \pmod{e_0},\\
e_1\cdot e_3 \equiv e_1+ e_3 \pmod{e_0}.
\end{align*}
Therefore, we obtain 
\begin{align} 
e_1 \cdot e_2&= e_1 + e_2+\widetilde{\beta_1}\,e_0, \label{e12}\\ 
e_2 \cdot e_3&= e_2 + e_3+\widetilde{\beta_2}\,e_0, \label{e23}\\ 
e_1 \cdot e_3&= e_1 + e_3+\widetilde{\beta_3}\,e_0, \label{e13}
\end{align}
for some $\widetilde{\beta_j}\ (j=1,2,3)$. 
Let ${\cal A}$ be a algebra determined by 
\eqref{ee1}, \eqref{e12}, \eqref{e23}, \eqref{e13}.
Here, we set 
\begin{align*}
\alpha_j&=\frac14 (1- \widetilde{\alpha_j}) \quad (j=1,2,3),\\
 c_1&=c_2=c_3=c_4=\frac1{4\gamma_3}(1+ \widetilde{\beta_3}),\\
\beta_1&=-\frac12 (1+\widetilde{\beta_1})- \frac{\gamma_1}{2\gamma_3}(1+ 
\widetilde{\beta_3}),\\ 
\beta_2&=-\frac12 (1+\widetilde{\beta_2})- \frac{\gamma_2}{2\gamma_3}(1+ 
\widetilde{\beta_3}), 
\end{align*}
and we determine a algebra ${\cal A}_3(c)$ associated with generalized 
Darboux-Halphen-Brioshi system \eqref{gDHB}. Then ${\cal A}_3(c)$ is 
equivalent with ${\cal A}$. Therefore ${\cal D}_3$ is a generalized 
Darboux-Halphen-Brioshi system. 

The converse is proved by a straightforward calculation. 
\end{proof}


\begin{thebibliography}{A}

\bibitem [CO]{CO}
Clarkson, P. A. and Olver, P. J.: 
{ Symmetry and the Chazy equation},
{\it J. Differential Equations} {\bf 124} (1996),  225--246. 

\bibitem [C]{C}
{Chazy, J.}: { Sur les \'equations  diff\'erentielles 
du trousi\`eme ordre et d'ordre sup\'erieur dont l'int\'egrale 
g\'en\'erale a ses points critiques fixes}, 
{\it Acta Math.} {\bf 34} (1911), 317--385.

\bibitem [D]{Dub}
Dubrovin, B.: { Painlev\'e trascendents in two-dimensional topological 
field theory}, 
{\it Proceedings of the 1996 Cargese sumer school "The Painlev\'e 
property: one century later"}, 
287--412 (1999).
 

\bibitem[H1]{H1}
{Halphen, G.}: { Sur une syst\'em d'\'equations diff\'erentielles},
{\it C. R. Acad. Sci., Paris} {\bf 92} (1881), 1101--1103.

\bibitem [H2]{H2}
{ G.~Halphen}, { Sur certains syst\`eme d'\'equations diff\'erentielles},
{\it C. R. Acad. Sci., Paris} {\bf 92} (1881), 1404--1406.

\bibitem [HM]{HM}
Harnad, J. and McKay, J.:
 Modular solutions to equations of generalized Halphen type, 
{\it R. Soc. Lond. Proc. Ser. A Math. Phys. Eng. Sci.} {\bf 456} (2000), 
261--294. 

\bibitem [KS]{KS}
Kinyon M. K. and Sagle A. A.:
{ Quadratic Dynamical Systems and Algebras}, 
{\it J. Differential Equations} {\bf 117} (1995)  67--126

\bibitem [M]{M} 
{Markus, L.}: { Quadratic differential equations and non-associative 
algebras}, 
{\it Ann.~Math. Stud.}, {\bf 45}, Princeton Univ.~Press, 1960, 185--213.

\bibitem [O]{O} 
{Ohyama, Y.}: 
Differential Relations of Theta Functions, {\it Osaka J. Math.} {\bf 32} 
431--450 (1995) 

\bibitem [O1]{O1} 
{Ohyama, Y.}: 
{ Nonlinear systems related to second order linear equations}, 
{\it Osaka J. Math.} {\bf 33} (1996), 927--949. 

\bibitem [O2]{O2} 
{Ohyama, Y.}: 
{ Hypergeometric functions and 
non-associative algebras}, 
{\it CRM Proceedings and Lecture Notes} {\bf 30} (2001) 173--184.

\bibitem [O3]{O3}
{Ohyama, Y.}: 
{ Differential equations for modular forms with level three}, 
{\it Funkcial. Ekvac.} {\bf 44} (2001) 377--389. 


\end{thebibliography}
\end{document}